\documentclass[12pt,leqno,fleqn]{amsart}  
\usepackage{amsmath,amstext,amsthm,amssymb,amsxtra} 

\usepackage{float}
\usepackage{amsrefs}
\usepackage{hyperref} 
\hypersetup{
    colorlinks=true,       
    linkcolor=blue,          
    citecolor=magenta,        
    filecolor=magenta,      
    urlcolor=cyan,          
}

\usepackage{pgf,pgfarrows} 
\usepackage{tikz}
\usetikzlibrary[decorations.pathreplacing]

\allowdisplaybreaks
\swapnumbers

\allowdisplaybreaks
\swapnumbers

\theoremstyle{plain} 
\newtheorem{lemma}[equation]{Lemma} 
 
\newtheorem{theorem}[equation]{Theorem}

\theoremstyle{definition}
\newtheorem{definition}[equation]{Definition} 

\theoremstyle{remark}

\numberwithin{equation}{section}

\def\norm#1.#2.{\lVert#1\rVert_{#2}}
\def\Norm#1.#2.{\bigl\lVert#1\bigr\rVert_{#2}}
\def\NOrm#1.#2.{\Bigl\lVert#1\Bigr\rVert_{#2}}
\def\NORm#1.#2.{\biggl\lVert#1\biggr\rVert_{#2}}
\def\NORM#1.#2.{\Biggl\lVert#1\Biggr\rVert_{#2}}

\def\ip#1,#2,{\langle #1,#2\rangle}
\def\Ip#1,#2,{\langle#1,#2\rangle}
\def\IP#1,#2,{\langle#1,#2\rangle}

\def\XXint#1#2#3{{\setbox0=\hbox{$#1{#2#3}{\int}$}
     \vcenter{\hbox{$#2#3$}}\kern-.5\wd0}}

\newcommand{\nint}{ \displaystyle\int}

\newcommand{\norme}[1]{ \left\| #1 \right\|}

\newcommand{\ral}{\mathbb{R}}

\newcommand{\nat}{\mathbb{N}}

\newcommand{\zat}{\mathbb{Z}}

\newcommand{\sop}[1]{{\mathbf{#1}}}
\newcommand{\vop}[1]{\mathcal{#1}_r}
\newcommand{\vecc}[1]{{\mathbf{#1}}}
\newcommand{\kapdep}{\kappa^4}

\title{A Note On $\ell^r$-Valued Calder\'{o}n-Zygmund Operators}
\author[J. Scurry]{James Scurry}
\address{ School of Mathematics, Georgia Institute of Technology, Atlanta GA 30332, USA}
\email {jscurry3@math.gatech.edu}
\begin{document}
\maketitle
\begin{abstract}
We consider $\ell^r$ extensions of Calder\'{o}n-Zygmund operators on weighted $L^p$ spaces. Our interest is in generalizing the scalar estimates for these operators (see \cite{hhhh}, \cite{npv}, \cite{handlandp}, and \cite{lerner2}) and the vector-valued theory considered by \cite{perez} and \cite{lernerlps}. In particular, we use multiple applications of Lerner's inequality to show that if $T$ is an $L^2(\ral^n)$ bounded Calder\'{o}n-Zygmund operator, its $\ell^r$ extension $\vecc{T}$ satisfies $\norme{\vecc{T}}_{L^p_{\ell^r}(w) \rightarrow L^p_{\ell^r}(w)} \lesssim [w]_{A_p}^{\max \{ 1 , \frac{1}{p-1} \}}$. Related results in general Banach spaces were studied in \cite{ht}.
\end{abstract}
\noindent
\begin{section}{Introduction}
\indent We intend to study $\ell^r$ extensions of Calder\'{o}n-Zygmund operators on weighted spaces $L^p(w)$ with $1 < p,r < \infty$ and $w \in A_p$. Our goal is to give a quantitative estimate of these operators' norm in terms of a given weight's $A_p$ characteristic. The scalar version of our problem has been given a great deal of attention. In this context the sharp dependence can be extrapolated from the case $p=2$ which gives a linear estimate, i.e. if $T$ is an $L^2(\ral^n)$ bounded Calder\'{o}n-Zygmund operator and $w \in A_2$,
\begin{eqnarray}
\norme{T}_{L^2(w) \rightarrow L^2(w)} &\lesssim& [w]_{A_2} \label{e.a2};
\end{eqnarray}
further, $(\ref{e.a2})$ is referred to as the $A_2$ Theorem. The authors of \cite{npv} reduced the proof of $(\ref{e.a2})$ to estimating Sawyer-type testing conditions: for $w \in A_2$,
\begin{eqnarray}
\norme{T}_{L^2(w) \rightarrow L^2(w)} &\lesssim& [w]_{A_2} + \norme{T}_{L^2(w) \rightarrow L^{2,\infty}(w)} + \norme{T^{\ast}}_{L^2(w^{-1}) \rightarrow L^{2,\infty}(w^{-1})} \label{e.a2222}.
\end{eqnarray}
Using probabilistic techniques, Hyt\"{o}nen first proved $(\ref{e.a2})$ in all generality by demonstrating the weak-type norms in $(\ref{e.a2222})$ satisfy a linear bound. 
Several subsequent proofs of $(\ref{e.a2})$ have also appeared, some of which appeal to averaging techniques (\cite{handpandv}, \cite{lerner2}) and others avoiding this altogether (\cite{handlandp}, \cite{lerner22}). \\
\indent In the vector-valued setting, several different types of operators have been considered. In \cite{perez}, 
the authors show the dyadic square function $S$ and vector-valued maximal operator $\vecc{M}_r$ with exponent $r$ satisfy:
\begin{eqnarray*}
\norme{S}_{L^p(w) \rightarrow L^p(w)} &\lesssim& [w]_{A_p}^{\max \{ \frac{1}{2}, \frac{1}{p-1} \} } \\
\norme{\vecc{M}_{r}}_{L^p_{\ell^r}(w) \rightarrow L^p(w)} &\lesssim& [w]_{A_p}^{\max \{ \frac{1}{r}, \frac{1}{p-1} \} }
\end{eqnarray*}
where $1 < p < \infty$ and $w \in A_p$. Using similar methods, \cite{lernerlps} gives sharp bounds for the intrinsic square function $G_{\alpha}$ on weighted $L^p(w)$ spaces, resolving a well-known conjecture. We aim to generalize the forgoing types of results to vector-valued extensions of an $L^2(\ral^n)$ bounded Calder\'{o}n-Zygmund operator, and the main theorem of this paper can be formulated as the following:
\begin{theorem}\label{t.main1}
Given a Calder\'{o}n-Zygmund operator $T$ on $\ral^n$, for $1 < r < \infty$ we denote by $\vecc{T}$ the $\ell^r$ extension of $T$, i.e. $\vecc{T}(\vecc{f}) = \{ T(f_j)(x) \}$ and 
\begin{eqnarray*}
\vop{T}(\vecc{f})(x) &=& \left( \displaystyle\sum_{j=1}^{\infty} |T(f_j)(x)|^r \right)^{\frac{1}{r}}
\end{eqnarray*} 
for $\vecc{f} = \{ f_j \}$ with $f_j \in \mathcal{S}(\ral^n)$.
Let $1 < p < \infty$ and $w \in A_p$. Given a Calder\'{o}n-Zygmund operator $T$ we have the following bounds:
\begin{eqnarray}
\norme{\vop{T}}_{L_{\ell^r}^{p}(w) \rightarrow L^{p}(w)} &\lesssim& [w]_{A_p}^{\max \{ 1, \frac{1}{p-1} \} } \label{m.sttype}.
\end{eqnarray}
\end{theorem}
\noindent
Unexpectedly, the strong type operator norm of $\vop{T}$ does not depend on $r$, indicating that scalar and vector-valued Calder\'{o}n-Zygmund operators can be equally singular. Additionally, the paper \cite{ht} considers more general Banach valued Calder\'{o}n-Zygmund operators and achieves our Theorem \ref{t.main1} as a corollary using different proof methods.   \\
\indent In the scalar case, the proof strategy is to reduce the study of $T$ to simpler operators, typically Haar-shift operators of a fixed complexity. We follow this tract, reducing the study of a given $T$ to consideration of vector-valued Haar-shift operators of a fixed complexity $\kappa$: indeed, we show it will be enough to prove the following theorem:
\begin{theorem} \label{t.main2}
Given a vector-valued Haar-shift operator $\vop{S}$ of complexity $\kappa$, we have
\begin{eqnarray*}
\norme{\vop{S}}_{L^p_{\ell^r}(w) \rightarrow L^p(w)} &\lesssim& \kapdep [w]_{A_p}^{\max \{ 1, \frac{1}{p-1} \}}.
\end{eqnarray*}
\end{theorem} 
\noindent
The chief difficulty in proving Theorem \ref{t.main2} will be maintaining a polynomial dependence on $\kappa$. As in \cite{lerner2}, \cite{lerner22}, and \cite{landh} we rely heavily on the application of Lerner's decomposition theorem. Specifically, we apply this inequality multiple times; the first application being component-wise to permit a decomposition of the resulting vector-valued operator analogous to the scalar decomposition considered in \cite{landh}. Then we follow \cite{lerner22} and apply Lerner's formula again to reduce our problem to vector-valued operators of complexity $1$; a third application of the theorem reduces our problem to the scalar case and completes the proof. \\
\indent The outline of the paper is as follows. In \ref{s.prelim} we introduce definitions and the main theorems of this paper; \ref{s.lt} lists several lemmas and theorems which will be used in our proofs. Subsequent sections refer to the proofs of specific theorems, beginning with arguments for our Lebesgue estimates and continuing with proofs of Theorem \ref{t.main2} and Theorem \ref{t.main1}.  
\subsubsection{Acknowledgment} The author would like to thank Dr. Michael Lacey for introducing the problem as well as crucial discussions and numerous suggestions. Further, the author would also like to thank Dr. Brett Wick for discussions concerning this paper, suggestions, and time.   
\begin{subsection}{Preliminaries} \label{s.prelim} 
In this section we fix notation and introduce our theorems. Let $1 < p, r < \infty$ and  $w \in A_p$ weight with $\kappa \in \nat$.  
\begin{definition}
For $u \in \{0,3^{-1} \}^n$ we denote by $\mathcal{D}^{u}$ the dyadic grid defined by
\begin{eqnarray*}
\mathcal{D}^{u} &=& \{ 2^{-k} ([0,1)^n + m +(-1)^k u) : k \in \zat, m \in \zat^n \}
\end{eqnarray*}
and note that this defines a collection of $2^n$ dyadic grids on $\ral^n$. In the special case $u = \vecc{0}$, we let $\mathcal{D}^{u} = \mathcal{D}$.
\end{definition}
\begin{definition}
We refer to a collection of cubes $\mathcal{Q} = \{ Q^k_j \}$ as sparse if 
\begin{itemize}
\item[i.] for fixed $j$, the $Q^k_j \cap Q^l_j = \emptyset$ \\
\item[ii.] for $Q^k_j \in \mathcal Q$ taking $D(Q^k_j) = Q^k_j \backslash \displaystyle\bigcup_{\substack{Q^m_l \subset Q^k_j \\ Q^m_l \in \mathcal Q}} Q^m_l$ we have
$|D(Q^k_j) \cap Q^k_j | \leq 2^{-1} |Q^k_j|$. 
\end{itemize}
\end{definition}

\begin{definition}
We call an operator $T$ a Calder\'{o}n-Zygmund operator in $\ral^n$ if $T$ is an $L^2(\ral^n)$ bounded integral operator with a kernel $K$ satisfying:
\begin{itemize}
\item[i.] $|K(x,y)| \lesssim \frac{1}{|x-y|^n}$ for $x,y \in \ral^n$ such that $x \neq y$
\item[ii.] $|K(x,y) - K(x^{\prime},y)|+ |K(y,x) - K(y,x^{\prime})| \lesssim \frac{|x-x^{\prime}|^{\alpha}}{|x-y|^{n+\alpha}}$ with $|x-x^{\prime}| < \frac{|x-y|}{2}$.
\end{itemize} 
\end{definition}
\begin{definition}
Let ${\bf{S}} = \{ S^j \}_{j=1}^{\infty}$ be a collection of generalized Haar shift operators of complexity $\kappa$ such that
$S^jf(x) = \displaystyle\sum_{I \in \mathcal D} \langle f, k^j_I \rangle h^j_I(x) = \displaystyle\sum_{I \in \mathcal D} S^j_If(x)$ for $f \in L^1_{\rm{loc}}(\ral^n)$. Take
\begin{eqnarray*}
{\vop{S}}{\mathbf{f}}(x) &=& \left( \displaystyle\sum_{j=1}^{\infty} |S^jf_j(x)|^r \right)^{\frac{1}{r}}.
\end{eqnarray*}
for $\vecc{f} = \{ f_j \}$ with $f_j \in L^1_{\rm{loc}}(\ral^n)$. We call $\vop{S}$ a vector-valued Haar-shift operator of complexity $\kappa$.
\end{definition}
\begin{definition}
We define an operator $\mathcal{P}_r$ as follows. For each $j$ let $\mathcal{Q}_j$ be a sparse collection of dyadic cubes from the same dyadic system. For $\vecc{f} = \{ f_j \}_{j=1}^{\infty}$, define
\begin{eqnarray*}
P^j(f_j)(x) &=& \displaystyle\sum_{Q \in \mathcal Q_j} \mathbb{E}_Q(f_j) \mathbf{1}_{E_j(Q)}(x)
\end{eqnarray*}
where for each $Q$, $E_j(Q)$ is a union of subcubes of $Q$ satisfying $ 2^{-\kappa } |Q| \leq |E_j(Q)|$ and take
\begin{eqnarray*}
\mathcal{P}_r (\vecc{f})(x) &=& \left( \displaystyle\sum_{j=1}^{\infty} |P^jf_j(x)|^r \right)^{\frac{1}{r}} .
\end{eqnarray*}
\noindent We refer to operators of the above type as positive vector-valued Haar-shift operators.
\end{definition}
\begin{definition}
For given $f \in L^1_{\rm{loc}}(\ral^n)$, $0 < \lambda < 1$, and $Q$ we have
\begin{eqnarray*}
\omega_{\lambda}(f; Q) &=& \displaystyle\inf_{c \in \ral} ((f - c) \mathbf{1}_{Q})^{\ast}(\lambda|Q|) \\
M^{\sharp}_{\lambda, Q}f(x) &=& \displaystyle\sup_{I \subset Q} \mathbf{1}_{Q}(x) \omega_{\lambda}(f, I)
\end{eqnarray*}
where for $g \in L^1_{\rm{loc}}(\ral^n)$, $g^{\ast}$ represents the symmetric non-increasing rearrangement. 
\end{definition}
\noindent
Now we list the main theorems of this paper:
\begin{theorem} \label{op.wktype} 
The operator $\vop{S}(\cdot)$ satisfies $\norme{\vop{S}}_{L^1_{\ell^r} \rightarrow L^{1,\infty}} \lesssim \kappa^{1+\frac{1}{r}}$.
\end{theorem}
\begin{theorem} \label{t.2pt}
For $\mathcal{P}_r$ as above, the following inequalities hold for Lebesgue measure:
\begin{eqnarray*}
\norme{\vop{P}}_{L^p_{\ell^r} \rightarrow L^p} &\lesssim& \kappa^2 \kappa^{ \max \{  r , r^{\prime} \} }.
\end{eqnarray*}
\end{theorem}
\begin{theorem}\label{t.3pt} 
With $w$ and $p$ as above we have
\begin{eqnarray}
 \norme{\vop{S}}_{L^p_{\ell^r}(w) \rightarrow L^{p,\infty}(w)} &\lesssim& 2^{\kappa} [w]_{A_p} \label{s.wktype} \\
\norme{\vop{S}}_{L^p_{\ell^r}(w) \rightarrow L^p(w)} &\lesssim& \kapdep [w]_{A_p}^{ \max \{ 1 , \frac{1}{p-1} \}} \label{s.strtype}.
\end{eqnarray}
\end{theorem}
\begin{theorem}\label{t.main4}
Let $T$ be an $L^2(\ral^n)$ bounded Calder\'{o}n-Zygmund operator and $w \in A_p$ with $1 < p < \infty$. For $1 < r < \infty$, 
\begin{eqnarray*}
\norme{\vop{T}}_{L^p_{\ell^r}(w) \rightarrow L^p(w)} &\lesssim& [w]_{A_p}^{\max \{ 1, \frac{1}{p-1} \}} .
\end{eqnarray*}
\end{theorem}

\end{subsection}
\begin{subsection}{Technical Lemmas and Theorems} \label{s.lt}
We begin by stating some technical Lemmas and Theorems which will be used to initiate our proofs. 

\begin{theorem}[Lerner] \label{t.lerner}
Let $f \in L^1_{\rm{loc}}(\ral^n)$ and let $Q$ be a fixed cube. Then there exists a collection of dyadic cubes $\{ Q^k_j \}_{j,k \in \nat}$ such that 
\begin{itemize}
\item[i.] for each $k,j \in \nat$, we have $Q^k_j \subset Q$ 
\item[ii.] for almost every $x \in Q$, 
\begin{eqnarray*}
| f(x) - m_{f}(Q) | &\leq& 4 M^{\sharp}_{2^{-n-2};Q}f(x) + 4 \displaystyle\sum_k \displaystyle\sum_{j} \omega_{2^{-n-2}} (f;Q^k_j) \mathbf{1}_{Q^k_j}(x)
\end{eqnarray*}
\item[iii.] for fixed $k$, $Q^k_j \cap Q^k_i = \emptyset $ for $i \neq j$ 
\item[iv.] letting $\Omega_k = \displaystyle\bigcup_{j} Q^k_j$, we have $| \Omega_k \cap Q^k_j| \leq 2^{-1} |Q^k_j|$ and $\Omega_{k+1} \subset \Omega_{k}$. 
\end{itemize}  
\end{theorem}

\begin{lemma}[Lemma 3.1, \cite{perez}] \label{lemma.perez} Given a measurable function $f$ and $Q \in \mathcal D$, then for $0 < \lambda < 1$ and $ 0 < p < \infty$ we have
\begin{eqnarray*}
(f \mathbf{1}_{Q})^{\ast}(\lambda | Q |) &\leq& \frac{\norme{f}_{L^{p,\infty}(Q, |Q|^{-1} dx)}}{\lambda^{\frac{1}{p}}}.
\end{eqnarray*}
\end{lemma}
\begin{lemma} \label{l.max}
If $Q \in \mathcal D$ then
\begin{eqnarray*}
\omega_{\lambda}({\vop{S}}{\mathbf{f}};Q) &\lesssim& \kappa^{1+\frac{1}{r}} 2^{\kappa} \mathbb{E}_{Q^{(\kappa)}}( \| \vecc{f} \|_{\ell^r})  . 
\end{eqnarray*}
\end{lemma}
\begin{lemma}\label{cz} [Proposition 2.3, \cite{lerner22}]
Let $T$ be a Calder\'{o}n-Zygmund operator and $Q \subset \ral^n$ a cube. If $ 1 < p < \infty$ and $w \in A_p$ then  for $f \in L^p(w)$ 
\begin{eqnarray*}
\omega_{\lambda}(Tf;Q) &\lesssim& \displaystyle\sum_{m=0}^{\infty} \frac{1}{2^{m \delta}} \left( \frac{1}{|2^mQ|} \nint_{2^mQ} |f(y)| dy \right)
\end{eqnarray*}
\end{lemma}
\begin{lemma}[Lemma 2.4, \cite{landh}] \label{l.landh} 
If $S$ is a generalized Haar-shift operator of complexity $\kappa$ then we have
\begin{eqnarray*}
\omega_{\lambda}(Sf;Q)(\lambda |Q|) &\lesssim& \frac{\kappa \mathbb{E}_{Q} |f|}{\lambda}   + \frac{1}{\lambda} \displaystyle\sum_{j=1}^{\kappa} \mathbb{E}_{Q^{(j)}}|f|.
\end{eqnarray*}
\end{lemma}

\begin{theorem}[Theorem 1.12, \cite{lerner2}] \label{t.shvmax}
Let $1 < q, p < \infty$, $0 < \lambda < 1$, and assume that $f$ and $g$ are functions satisfying the following: for any cube $Q$ we have
\begin{eqnarray*}
\omega_{\lambda}(|g|^q; Q) &\lesssim& \left( \frac{ \mathbb{E}_{Q} |f|}{\lambda} \right)^q
\end{eqnarray*} 
for some constant independent of $Q$.
Then we have
\begin{eqnarray*}
\norme{g}_{L^p}(w) &\lesssim& [w]_{A_p}^{\max \{\frac{1}{q}, \frac{1}{p-1} \}} \norme{f}_{L^p(w)} .
\end{eqnarray*}
\end{theorem}
\begin{theorem}[Theorem 1.12, \cite{perez}] \label{t.maxv}
For $ 1 < r , p < \infty$ and $w \in A_p$ we have the following bound:
\begin{eqnarray*}
\norme{\vecc{M}_{r}}_{L^p_{\ell^r}(w) \rightarrow L^p(w)} &\lesssim& [w]_{A_p}^{\max \{ \frac{1}{r}, \frac{1}{p-1} \}}. 
\end{eqnarray*}
\end{theorem}
\end{subsection}
\end{section}
\begin{section}{The Lebesgue Estimates}
\begin{subsection}{Proof of Theorem \ref{op.wktype}}
We will perform a Calderon-Zygmund decomposition. Fix $\lambda > 0$ and let $\{ Q_j \}_{j=1}^{\infty}$ be the maximal dyadic cubes such that $\frac{1}{|Q_j|} \nint_{Q_j} \norme{\vecc{f}}_{\ell^r} dx \geq \lambda$. For each $j$ define $\vecc{b}^j$ by
\begin{eqnarray*}
b^j_k(x) &=& \left( f_k -  \frac{1}{|Q_j|} \nint_{Q_j} f_k \right) \mathbf{1}_{Q_j}(x).
\end{eqnarray*}
and let $\vecc{b} = \displaystyle\sum_{j=1}^{\infty} \vecc{b}^j$. Further, we let $\vecc{g} = \vecc{f} - \vecc{b}$. Then we have the following:
\begin{itemize}
\item[(i)] $\norme{\vecc{g}}_{L^1_{\ell^r}} \lesssim \norme{\vecc{f}}_{L^1_{\ell^r}}$
\item[(ii)] for each $j$, ${\rm{supp}} \ b^j_k \subset Q_j$ all $k \in \nat$ 
\item[(iii)] $\displaystyle\sum_{j=1}^{\infty} \norme{\vecc{b}^j}_{L^1_{\ell^r}} \lesssim \norme{\vecc{f}}_{L^1_{\ell^r}}$ 
\item[(iv)]  for almost all $x \in \ral$, $\norme{\vecc{g}}_{\ell^r} \lesssim \lambda \norme{\vecc{f}}_{\ell^r}$
\item[(v)] $\displaystyle\sum_{j=1}^{\infty} | Q_j | \lesssim \frac{\norme{\vecc{f}}_{L^1_{\ell^r}}}{\lambda}$.
\end{itemize} 
\noindent Notice 
\begin{eqnarray*}
| \{ x \in \ral^n : \vop{S}\vecc{f}(x) > \lambda \} | &\leq& \left| \left\{ x \in \ral^n : \vop{S}\vecc{g}(x) > \frac{\lambda}{2} \right\} \right| + 
\left| \left\{ x \in \ral^n : \vop{S}\vecc{b}(x) > \frac{\lambda}{2} \right\} \right|
\end{eqnarray*}
\noindent and consider by Chebyshev's inequality,
\begin{eqnarray*}
\left| \left\{ x \in \ral^n : \vop{S}\vecc{g}(x) > \frac{\lambda}{2} \right\} \right| &\leq& \frac{4}{\lambda^2} \nint_{\ral^n} \vop{S}\vecc{g}(x)^2 dx \\
&\lesssim& \frac{4}{\lambda^2} \nint_{\ral^n} \norme{\vecc{g}}_{L^2_{\ell^r}}.
\end{eqnarray*}
\noindent By properties (i) and (iv) from above, 
\begin{eqnarray*}
\nint_{\ral^n} \norme{\vecc{g}}_{\ell^r}^2 dx 
&\lesssim& \lambda \nint_{\ral^n} \norme{\vecc{f}}_{\ell^r} dx
\end{eqnarray*}
\noindent so that
\begin{eqnarray*}
\frac{4}{\lambda^2} \nint_{\ral^n} \norme{\vecc{g}}_{\ell^r} dx &\lesssim& \frac{4}{\lambda} \nint_{\ral^n} \norme{\vecc{g}}_{\ell^r} dx.
\end{eqnarray*}
\noindent On the other hand,
\begin{eqnarray*}
\vop{S}\vecc{b}(x) &\leq& \displaystyle\sum_{j=1}^{\infty} \vop{S}\vecc{b}^j(x).
\end{eqnarray*}
\noindent Further, for $(Q_j)^{(\kappa)} \subset I$, we have $\nint_{I} b^j_k(x) dx = 0$ so that 
$S^k_I(b^j_k)(x) =0$ for $(Q_j)^{(\kappa)} \subset I$. Hence, by standard computations 
\begin{eqnarray*}
\displaystyle\sum_{j=1}^{\infty} \vop{S}\vecc{b}^j(x)  &=&
\displaystyle\sum_{j=1}^{\infty} \left( \displaystyle\sum_{k=1}^{\infty} \left| S^k(b^j_k)(x) \right|^r \right)^{\frac{1}{r}} \\
&\leq& \displaystyle\sum_{j=1}^{\infty} \left( \displaystyle\sum_{k=1}^{\infty} \left| \displaystyle\sum_{I \subseteq Q_j} S^k_I(b^j_k)(x) \right|^r \right)^{\frac{1}{r}} + \kappa^{\frac{1}{r}} \displaystyle\sum_{j=1}^{\infty}    \displaystyle\sum_{Q_j \subset I \subset Q^{(\kappa)}_j} \left( \displaystyle\sum_{k=1}^{\infty} \left( \mathbb{E}_{I}|b^j_k| \right)^r \mathbf{1}_{I}(x)  \right)^{\frac{1}{r}} \\
&=& \displaystyle\sum_{j=1}^{\infty} \left( \displaystyle\sum_{k=1}^{\infty} \left| \displaystyle\sum_{I \subseteq Q_j} S^k_I(b^j_k)(x) \right|^r \right)^{\frac{1}{r}} + \kappa^{\frac{1}{r}} \displaystyle\sum_{Q_j \subset I \subset Q^{(\kappa)}_j} \displaystyle\sum_{j=1}^{\infty}     \left( \displaystyle\sum_{k=1}^{\infty} \left( \mathbb{E}_{I}|b^j_k| \right)^r \mathbf{1}_{I}(x) \right)^{\frac{1}{r}}
\end{eqnarray*}
\noindent Let
\begin{eqnarray*}
A &=& \displaystyle\sum_{j=1}^{\infty} \left( \displaystyle\sum_{k=1}^{\infty} \left| \displaystyle\sum_{I \subseteq Q_j} S^k_I(b^j_k)(x) \right|^r \right)^{\frac{1}{r}} \\
B &=& \kappa^{\frac{1}{r}} \displaystyle\sum_{Q_j \subset I \subset Q^{(\kappa)}_j} \displaystyle\sum_{j=1}^{\infty} \left( \displaystyle\sum_{k=1}^{\infty} \left( \mathbb{E}_{I}|b^j_k| \right)^r \mathbf{1}_{I}(x) \right)^{\frac{1}{r}} 
\end{eqnarray*}
\noindent so that
\begin{eqnarray*}
\left| \left\{ \vop{S} \vecc{b}(x) > \frac{\lambda}{2} \right\} \right| &\leq& 
\left| \left\{ A > \frac{\lambda}{4} \right\} \right| + \left| \left\{ B > \frac{\lambda}{4} \right\} \right|.
\end{eqnarray*}
\noindent Notice that $A$ is supported on $\displaystyle\cup Q_j$ so that
\begin{eqnarray*}
\left| \left\{ A > \frac{\lambda}{4} \right\} \right| \leq \displaystyle\sum_{j=1}^{\infty} | Q_j | \lesssim \frac{\norme{\vecc{f}}_{L^1_{\ell^r}}}{\lambda}
\end{eqnarray*} 
\noindent and using Chebyshev's inequality we have
\begin{eqnarray*}
\left| \left\{ B > \frac{\lambda}{4} \right\} \right| &\leq&
\frac{4}{ \lambda} \kappa^{\frac{1}{r}} \nint  \displaystyle\sum_{Q_j \subset I \subset Q^{(\kappa)}_j} \displaystyle\sum_{j=1}^{\infty} \left( \displaystyle\sum_{k=1}^{\infty} \left( \mathbb{E}_{I}|b^j_k| \right)^r \mathbf{1}_{I}(x)  \right)^{\frac{1}{r}} dx. \label{e.minkowskii} 
\end{eqnarray*}
\noindent Applying Minkowskii's integral inequality two the inner sum of expectations in (\ref{e.minkowskii}) yields
\begin{eqnarray*}
(\ref{e.minkowskii}) &\leq& \frac{4 \kappa^{\frac{1}{r}} }{ \lambda} \nint \displaystyle\sum_{Q_j \subset I \subset Q^{(\kappa)}_j}  \displaystyle\sum_{j=1}^{\infty}  \mathbb{E}_{I} (\norme{\vecc{b}^j}_{\ell^r}) \mathbf{1}_{I} (x)  dx \\
&\leq& \frac{4 \kappa^{\frac{1}{r}}}{\lambda} \displaystyle\sum_{Q_j \subset I \subset Q^{(\kappa)}_j} \displaystyle\sum_{j=1}^{\infty} \norme{\vecc{b}^j}_{L^1_{\ell^r}} \\
&\lesssim& \frac{4 \kappa^{\frac{1}{r}}}{\lambda} \displaystyle\sum_{Q_j \subset I \subset Q^{(\kappa)}} \norme{\vecc{f}}_{L^1_{\ell^r}} \\
&\leq& \frac{4 \kappa^{1+ \frac{1}{r}}}{\lambda} \norme{\vecc{f}}_{L^1_{\ell^r}}.
\end{eqnarray*}
\noindent Combining the above estimates gives $\norme{\vecc{S}}_{L^1_{\ell^r} \rightarrow L^{1,\infty}} \lesssim \kappa^{1+ \frac{1}{r}}$.
\end{subsection}
\begin{subsection}{Proof of Theorem \ref{t.2pt}}
Fix $\vecc{f} \in L^p_{\ell^r}$ and suppose first $p = r$. In this case we have
\begin{eqnarray*}
\nint_{\ral^n} \vop{P}(\vecc{f})(x)^p dx &=&
\nint_{\ral^n} \displaystyle\sum_{j=1}^{\infty} |P^j(f_j)(x)|^r dx \\
&\lesssim& \kappa^r \nint_{\ral^n} \displaystyle\sum_{j=1}^{\infty} |f_j(x)|^r dx \\
&=& \kappa^r \nint_{\ral^n} \norme{\vecc{f}}^p dx.
\end{eqnarray*}
Now by Theorem \ref{op.wktype} and the Marcinkiewicz Interpolation Theorem for vector-valued operators we have for $1< p \leq r$,
\begin{eqnarray*} 
\norme{\vop{P}}_{L^p_{\ell^r}\rightarrow L^p} &\lesssim& \kappa^{r+2}.
\end{eqnarray*}
For the range $1 < r < p$ we notice there is a vector $\vecc{h} \in L^{p^{\prime}}_{\ell^{r^{\prime}}}$ with $\norme{\vecc{h}}_{L^{p^{\prime}}_{\ell^{r^{\prime}}}} = 1$ such that
\begin{eqnarray*}
\nint_{\ral^n} \vop{P}(\vecc{f})(x)^p dx &=&
\nint_{\ral^n} \vop{P}(\vecc{f}) \cdot \vecc{h} dx \\
&\leq& \left( \nint_{\ral^n} \norme{\vecc{f}}_{\ell^r}^p dx \right)^{\frac{1}{p}} \left( \nint_{\ral^n} \vecc{U}(\vecc{h})(x)^{p^{\prime}} dx \right)^{\frac{1}{p^{\prime}}} 
\end{eqnarray*}
where $\vecc{U}$ represents a `dual' operator for $\vop{P}$, i.e. if $(P^j)^{\ast}$ is the dual for each $P^j$ then
\begin{eqnarray*}
\vecc{U}(\vecc{g})(x)&=& \left( \displaystyle\sum_{j=1}^{\infty} |(P^j)^{\ast}(g_j)(x)|^{r^{\prime}} \right)^{\frac{1}{r^{\prime}}}
\end{eqnarray*}
with $\vecc{g} = \{ g_j \}$ and $g_j \in L^1_{\rm{loc}}(\ral^n)$. Arguing as before with $\vecc{U}$ in place of $\vop{P}$, we see
\begin{eqnarray*}
\left( \nint_{\ral^n} \vecc{U}(\vecc{h})(x)^{p^{\prime}} dx \right)^{\frac{1}{p^{\prime}}}
&\lesssim& \kappa^{r^{\prime}+2}.
\end{eqnarray*}
Hence, we have
\begin{eqnarray*}
\norme{\vop{P}}_{L^p_{\ell^r} \rightarrow L^p} &\lesssim& \kappa^2 \max \{ \kappa^{r}, \kappa^{r^{\prime}} \}.
\end{eqnarray*}
\end{subsection}

\begin{subsection}{Proof of Theorem \ref{l.max}}
By the triangle inequality we have,
\begin{eqnarray*}
\left| \mathbf{1}_{Q}(x) \vop{S}(\sop{f})(x) - \mathbf{1}_{Q}(x) \vop{S}( \mathbf{1}_{(Q^{(\kappa)})^{\rm{c}}} \vecc{f})(x) \right| &\leq& \mathbf{1}_{Q}(x) \vop{S}(\vecc{f} \mathbf{1}_{Q^{(\kappa)}})(x).
\end{eqnarray*}
Notice, $\vop{S}(\vecc{f} \mathbf{1}_{(Q^{(\kappa)})^{\rm{c}}})(x) \mathbf{1}_{Q}(x)$ is constant on $Q$. Define
\begin{eqnarray*}
C(Q,{\mathbf{f}}, \kappa) =  C = \vop{S}\left( \mathbf{1}_{(Q^{(\kappa)})^{\rm{c}}} {\mathbf{f}} \right)(x) \mathbf{1}_Q(x).
\end{eqnarray*}
\noindent Now the above implies 
\begin{eqnarray*}
\omega_{\lambda}(\vop{S}(\vecc{f}); \lambda |Q|) &\leq& 
(\mathbf{1}_{Q} \vop{S}(\vecc{f} \mathbf{1}_{Q^{(\kappa)}}))^{\ast} \left( {\lambda |Q|} \right). \end{eqnarray*}
Applying Lemma \ref{lemma.perez} gives
 \begin{eqnarray*}
 (\mathbf{1}_{Q} \vop{S}(\vecc{f} \mathbf{1}_{Q^{(\kappa)}}))^{\ast} \left( \lambda |Q| \right) &\lesssim& 
 \norme{ \vop{S}(\vecc{f} \mathbf{1}_{Q^{(\kappa)}})}_{L^{1,\infty}(Q, |Q|^{-1} dx)} 
\end{eqnarray*}
and from the weak-(1,1) inequality for $\vop{S}$ we obtain
\begin{eqnarray*}
\norme{ \vop{S} \left(\vecc{f} \mathbf{1}_{Q^{(\kappa)}} \right)}_{L^{1,\infty}(Q, |Q|^{-1} dx)} &\lesssim&\kappa^{1+ \frac{1}{r}} 2^{\kappa} \mathbb{E}_{Q^{(\kappa)}} (\| \vecc{f} \|_{\ell^r})  \\
 \end{eqnarray*}
Thus, 
\begin{eqnarray*}
\omega_{\lambda}(\vop{S}(\vecc{f}); \lambda |Q|) &\lesssim& \kappa^{1+\frac{1}{r}} 2^{\kappa} \mathbb{E}_{Q^{(\kappa)}} \| \vecc{f} \|_{\ell^r}  
\end{eqnarray*}
\end{subsection}

\end{section}
\begin{section}{Proof of Theorem \ref{t.3pt}}
\begin{subsection}{Proof of \eqref{s.wktype}}
Let $\mathbf{f} \in L^p_{\ell^r}(w)$ such that $\norme{\mathbf{f}}_{\ell^r}$ has compact support. Recall, $m_{S^j(f_j)}(Q) \rightarrow 0$ as $\ell(Q) \rightarrow \infty$ and in particular, for all cubes $Q$ which are sufficiently large, we have the following point-wise bound:
\begin{eqnarray*}
\vop{S}\vecc{f}(x)^p &\lesssim& M^{\sharp}(\norme{\vecc{f}}_{\ell^r})(x)^p + 4 M^{\sharp}_{2^{-n-1};Q_N}(\vop{S}\vecc{f})(x) + 4 \displaystyle\sum_{I \in \mathcal K} \omega_{2^{-n-1}}(\vop{S}\vecc{f};I) \mathbf{1}_I(x),
\end{eqnarray*}
where $\mathcal{K}$ is a sparse collection of cubes.
By Lemma \ref{op.wktype}, 
\begin{eqnarray*}
M^{\sharp}_{2^{-n-1};Q}(\vop{S}\vecc{f})(x) &\lesssim& M(\norme{\vecc{f}}_{\ell^r})(x) \\
\omega_{2^{-n-1}}(\vop{S}\vecc{f};I) \mathbf{1}_I(x) &\lesssim& \mathbb{E}_{I^{(\kappa)}}(\norme{\vecc{f}}_{\ell^r}) \mathbf{1}_I(x)
\end{eqnarray*}
so that 
\begin{eqnarray}
\vop{S}\vecc{f}(x) &\lesssim& M(\norme{\vecc{f}}_{\ell^r})(x) +   \displaystyle\sum_{I \in \mathcal K} 2^{\kappa} \mathbb{E}_{I}(\norme{\vecc{f}}_{\ell^r}) \mathbf{1}_I(x) \\
&=& M(\norme{\vecc{f}}_{\ell^r})(x) + 2^{\kappa} S(\norme{\vecc{f}}_{\ell^r})(x) \label{e.lernpwse}.
\end{eqnarray}
As a result, we have
\begin{eqnarray*}
w \left( \{ x \in \ral^n : \vop{S}\vecc{f}(x) > \alpha \} \right) &\leq&
w \left( \{ x \in \ral^n : M(\norme{\vecc{f}}_{\ell^r})(x) \gtrsim \frac{\alpha}{2} \} \right) + \\
\qquad & & 
w \left( \{ x \in \ral^n : 2^{\kappa} S(\norme{\vecc{f}}_{\ell^r})(x) \gtrsim \frac{\alpha}{2}  \} \right).
\end{eqnarray*}
By Buckley's bound (\cite{buckley}),
\begin{eqnarray}
w \left( \{ x \in \ral^n : M (\norme{\vecc{f}}_{\ell^r})(x) \gtrsim \frac{\alpha}{2} \} \right)
&\lesssim& \frac{[w]_{A_p}}{\alpha^p} \nint_{\ral^n} \norme{\vecc{f}}^p w \label{e.buckk} 
\end{eqnarray}
and since $S$ is a Haar-shift operator of complexity $\kappa$, we have
\begin{eqnarray}
w \left( \{ x \in \ral^n : 2^{\kappa} S( \norme{\vecc{f}}_{\ell^r} )(x) \gtrsim \frac{\alpha}{2} \right)
&\lesssim& \frac{[w]_{A_p}^p 2^{\kappa p}}{\alpha^p} \nint_{\ral^n} \norme{\vecc{f}}^p w \label{e.lsut} ; 
\end{eqnarray}
combining $(\ref{e.buckk})$ and $(\ref{e.lsut})$ gives the weak-type bound.
\end{subsection}
\begin{subsection}{Proof of \eqref{s.strtype}}
Let $\mathbf{f} \in L^p_{\ell^r}(w)$ such that $\norme{\mathbf{f}}_{\ell^r}$ has compact support. By applying Lerner's inequality to each component of $\vop{S}$ on a sufficiently large cube $J$, we obtain the bound 
\begin{eqnarray*}
\vop{S}(\vecc{f})(x) &\lesssim&
\left( \displaystyle\sum_{j=1}^{\infty} M^{\sharp}_{\frac{1}{4}; J}(S^j(f_j))(x)^r 
\right)^{\frac{1}{r}} +  \left( \displaystyle\sum_{j=1}^{\infty}  \left(\displaystyle\sum_{Q \in \mathcal Q_j} \mathbf{1}_{Q}(x)
\omega_{2^{-n-2}} (S^j(f_j); Q)  \right)^r \right)^{\frac{1}{r}} 
\end{eqnarray*}
\noindent
where $\mathcal{Q}_j$ is the collection of cubes which results from applying Theorem \ref{t.lerner} to $S^j(f_j)$.
\noindent Using Lemma \ref{l.landh} as in \cite{landh} we obtain for each $j$,
\begin{eqnarray*}
M^{\sharp}_{\frac{1}{4}; J}(S^j(f_j))(x)^r &\lesssim& \kappa^r Mf_j(x)^r
\end{eqnarray*}
and 
\begin{eqnarray*}
\left( \displaystyle\sum_{j=1}^{\infty} M^{\sharp}_{\frac{1}{4}; J}(S^j(f_j))(x)^r 
\right)^{\frac{1}{r}} &\lesssim& \vecc{M}_r(\vecc{f})(x).
\end{eqnarray*}
Now we consider the function
\begin{eqnarray}
\left( \displaystyle\sum_{j=1}^{\infty}  \left(\displaystyle\sum_{Q \in \mathcal Q_j} \mathbf{1}_{Q}(x) 
\omega_{2^{-n-2}} (S^j(f_j); Q)  \right)^r \right)^{\frac{1}{r}}. \label{e.subr}
\end{eqnarray}
Applying Lemma \ref{l.landh} for each $j$ we obtain
\begin{eqnarray}
(\ref{e.subr}) &\lesssim& \left( \displaystyle\sum_{j=1}^{\infty}  \left(\displaystyle\sum_{Q \in \mathcal Q_j} \kappa \cdot \mathbf{1}_{Q}(x) \cdot \mathbb{E}_{Q} |f_j|   +  \displaystyle\sum_{i=1}^{\kappa} \mathbf{1}_{Q}(x) \cdot \mathbb{E}_{(Q)^{(i)}}|f_j|   \right)^r \right)^{\frac{1}{r}} \label{e.aplerner}
\end{eqnarray} 
For each $j$ and $0 \leq i \leq \kappa$ define $E(Q)^{i} = \displaystyle\bigcup_{ \substack{(I)^{(i)} = Q \\ I \in \mathcal Q_j}} I$ with the convention $E(Q)^{0} = Q$. Then we recall, if $\{ x_i \}_{i=1}^n$ is a non-negative sequence of numbers, for $0 < q < \infty$,
\begin{eqnarray}\label{bb.i}
\left( \displaystyle\sum_{i=1}^n x_i \right)^q &\leq& n^q \displaystyle\sum_{i=1}^n x_i^q;
\end{eqnarray} 
applying $(\ref{bb.i})$ we obtain
\begin{eqnarray*}
(\ref{e.aplerner})
&\lesssim& \kappa^2  \left( \displaystyle\sum_{j=1}^{\infty}  \left( \displaystyle\sum_{Q \in \mathcal Q_j}  \mathbf{1}_{E(Q)^0}(x) \cdot (\mathbb{E}_{Q} |f_j|) \right)^r \right)^{\frac{1}{r}} + \\
& & \qquad{} \kappa^2 \displaystyle\sum_{i=1}^{\kappa}  \left( \displaystyle\sum_{j=1}^{\infty} \left( \displaystyle\sum_{Q \in \mathcal Q_j} \mathbf{1}_{E(Q)^i}(x) \cdot (\mathbb{E}_{(Q)^{(i)}}|f_j|) \right)^r \right)^{\frac{1}{r}}. 
\end{eqnarray*}
For $0 \leq i \leq \kappa$ let
\begin{eqnarray*}
P^{j,i}(g)(x) &=& \displaystyle\sum_{Q \in \mathcal Q_j}  \mathbf{1}_{E(Q)^i}(x) \mathbb{E}_{Q}(g)  
\end{eqnarray*}
with $g \in L^1_{\rm{loc}}(\ral^n)$ and $\vecc{P}_i$ be defined by
\begin{eqnarray*}
\vecc{P}_i(\vecc{g})(x) &=& \left(\displaystyle\sum_{j=1}^{\infty} |P^{j,i}(g)(x)|^r \right)^{\frac{1}{r}}
\end{eqnarray*}
for $\vecc{g} = \{ g_j \}$ and $g_j \in L^1_{\rm{loc}}(\ral^n)$. Hence we have the following point-wise bound
\begin{eqnarray}
\vop{S}(\vecc{f})(x) &\lesssim&  \vecc{M}_r(\vecc{f})(x) + \kappa^2 \displaystyle\sum_{i=0}^{\kappa} \vecc{P}_i(\vecc{f})(x) \label{e.ptwse} . 
\end{eqnarray}
By Theorem \ref{t.maxv},
\begin{eqnarray*}
\nint_{\ral^n} \vecc{M}_r(\vecc{f})(x)^p w
&\lesssim& [w]_{A_p}^{\max \{ \frac{p}{r}, \frac{p}{p-1} \} } \norme{\vecc{f}}_{L^p_{\ell^r}(w)}^p.
\end{eqnarray*}
So from $(\ref{e.ptwse})$, 
\begin{eqnarray*}
\nint_{\ral^n} \vop{S}(\vecc{f})(x)^p w &\lesssim& 
 \nint_{\ral^n} \vecc{M}_r(\vecc{f})(x)^p w + \kappa^{2p} \displaystyle\sum_{i=0}^{\kappa} \nint_{\ral^n} \vecc{P}_i(\vecc{f})(x)^p w \\
&\lesssim& [w]_{A_p}^{\max \{ p, \frac{p}{r} \} } \nint_{\ral^n} \norme{\vecc{f}}_{\ell^r}^p w + \kappa^{2p} \displaystyle\sum_{i=0}^{\kappa} \nint_{\ral^n} \vecc{P}_i(\vecc{f})(x)^p w.
\end{eqnarray*}
From duality, there is a vector $\vecc{h} = \{ h_j \} \in L^{p^{\prime}}_{\ell^{r^{\prime}}}(w)$ such that
\begin{eqnarray*}
\left( \nint_{\ral^n} \vecc{P}_i(\vecc{f})(x)^{p} w \right)^{\frac{1}{p}} &=&
\nint_{\ral^n} \vecc{P}_i (\vecc{f})(x) \cdot \vecc{h} w \\ 
&\leq& \norme{\vecc{f}}_{L^p_{\ell^r}(w)} \left( \nint_{\ral^n} \vecc{U}_i(\vecc{h} w)(x)^{p^{\prime}} \sigma \right)^{\frac{1}{p^{\prime}}},
\end{eqnarray*}
where for each $i$
\begin{eqnarray*}
\vecc{U}_i(\vecc{g})(x) &=& \left( \displaystyle\sum_{j=1}^{\infty} |(P^{j,i})^{\ast}(g_j)(x)|^{r^{\prime}} \right)^{\frac{1}{r^{\prime}}}
\end{eqnarray*}
with $\vecc{g} = \{ g_j \}$ and $g_j \in L^1_{\rm{loc}}(\ral^n)$.
We apply Lerner's Theorem in each component of $\vecc{U}_i$ to obtain the bound
\begin{eqnarray*}
\vecc{U}_i(\vecc{h}w)(x) &\lesssim& \vecc{M}_{r^{\prime}}(\vecc{h}w)(x) + \kappa \left( \displaystyle\sum_{j=1}^{\infty} |L^j(h_jw) |^{r^{\prime}} \right)^{\frac{1}{r^{\prime}}} \\
&=& \vecc{M}_{r^{\prime}}(\vecc{h}w)(x) + \kappa \mathcal{L}(\vecc{h}w)(x)
\end{eqnarray*} 
where for each $j$,
\begin{eqnarray*}
L^j(h_jw)(x) &=& \displaystyle\sum_{I \in \mathcal L_j} \mathbb{E}_I(h_jw) \mathbf{1}_I(x).
\end{eqnarray*}
Notice $\mathcal{L}$ is a vector-valued Haar-shift operator of complexity $1$ which is $L^2(\ral^n)$ bounded; hence, by Theorem \ref{l.max},
\begin{eqnarray*}
\omega_{\lambda}(\mathcal{L}(\vecc{h}w))(\lambda |Q|) &\lesssim& \mathbb{E}_{Q}(\norme{\vecc{h}}_{\ell^{r^{\prime}}}w)
\end{eqnarray*}
so that from another application of Lerner's Theorem we obtain a sparse collection of cubes $\mathcal {K}$,
\begin{eqnarray*}
\mathcal{L}(\vecc{h}w)(x) &\lesssim& M(\norme{\vecc{h}}_{\ell^{r^{\prime}}}w)(x) +\displaystyle\sum_{I \in \mathcal K} \mathbb{E}_{I}(\norme{\vecc{h}}_{\ell^{r^{\prime}}}w) \mathbf{1}_{I}(x).
\end{eqnarray*}
Hence for each $i$ we have 
\begin{eqnarray*}
\nint_{\ral^n} \vecc{U}_i(\vecc{h} w)(x)^{p^{\prime}} \sigma &\lesssim& \nint_{\ral^n} \vecc{M}_{r^{\prime}}(\vecc{h} w)(x)^{p^{\prime}} \sigma + \kappa^{p^{\prime}} \nint_{\ral^n} M(\norme{\vecc{h}}_{\ell^{r^{\prime}}}w)^{p^{\prime}} \sigma + \kappa^{p^{\prime}} \nint_{\ral^n} \mathcal{L}(\norme{\vecc{h}}_{\ell^r} w)(x)^{p^{\prime}} \sigma \\
&\lesssim& \kappa^{p^{\prime}} [w]_{A_p}^{\max \{ p^{\prime}, \frac{p^{\prime}}{p-1} \} } \nint_{\ral^n} \norme{\vecc{h}}_{\ell^{r^{\prime}}}^{p^{\prime}} \sigma \\
&\lesssim& \kappa^{p^{\prime}} [w]_{A_p}^{\max \{ p^{\prime}, \frac{p^{\prime}}{p-1} \} } .
\end{eqnarray*}
Now,
\begin{eqnarray*}
\norme{\vop{S}}_{L^p_{\ell^r}(w) \rightarrow L^{p}(w)} &\lesssim& \kapdep[w]_{A_p}^{\max \{ 1 , \frac{1}{p-1} \} }. 
\end{eqnarray*}
giving the result.
\end{subsection}
\begin{subsection}{An Example}
The bound for $\vop{S}$ is sharp by the scalar bound, but here we give an explicit example to show the bound is sharp. For each $j$ let $I_j = [0,2^{-j})$ and define
\begin{eqnarray*}
S(f)(x) &=& \displaystyle\sum_{j=1}^{\infty} \mathbb{E}_{I_j}(f) \mathbf{1}_{I_j}(x).
\end{eqnarray*}
Let $w(x) = |x|^{(\delta-1)(p-1)}$ and $f(x) = |x|^{\delta - 1}\mathbf{1}_{[0,1)}(x)$. Then
\begin{eqnarray*}
\norme{f}_{L^p(w)}^p &=& \nint_{[0,1)} |x|^{\delta-1} dx \\
&=& \frac{1}{\delta}.
\end{eqnarray*}
On the other hand
\begin{eqnarray*}
\norme{S(f)}_{L^p(w)}^p &=&
\nint_{[0,1)} \left( \displaystyle\sum_{j=1}^{\infty} \mathbb{E}_{I_j}(|x|^{\delta-1}) \mathbf{1}_{I_j}(x) \right)^{p} |x|^{(1-\delta)(p-1)} dx \\
&=& 
\displaystyle\sum_{k = 0}^{\infty} \nint_{[2^{-k-1},2^{-k})} \left( \displaystyle\sum_{j=0}^{\infty} \mathbb{E}_{I_j}(|x|^{\delta-1}) \mathbf{1}_{I_j}(x) \right)^p |x|^{(1-\delta)(p-1)} dx \\
&\sim& \displaystyle\sum_{k = 0}^{\infty} \nint_{[2^{-k-1},2^{-k})} \delta^{-p} |x|^{(\delta-1)p} |x|^{(1-\delta)(p-1)} dx \\
&=& \nint_{[0,1)} \delta^{-p} |x|^{\delta - 1} dx\\ 
&=& \delta^{-p-1} .
\end{eqnarray*}
Hence,
\begin{eqnarray*}
[w]_{A_p}^{\frac{1}{p-1}} &\sim& \delta^{-1} \\
&\lesssim& \frac{\norme{S(f)}_{L^p(w)}}{\norme{f}_{L^p(w)}}.
\end{eqnarray*}
As a consequence,
\begin{eqnarray*}
\left( \nint S(f \sigma)(x)^p w \right)^{\frac{1}{p}} &=& \displaystyle\sup_{\substack{ h \in L^{p^{\prime}}(w) \\ \norme{h}_{L^{p^{\prime}}(w)} = 1}} \nint S(f \sigma)(x) h(x) w \\ &=& \displaystyle\sup_{\substack{ h \in L^{p^{\prime}}(w) \\ \norme{h}_{L^{p^{\prime}}(w)} = 1}} \nint f(x) S^{\ast}(hw)(x) \sigma \\
&=& \displaystyle\sup_{\substack{ h \in L^{p^{\prime}}(w) \\ \norme{h}_{L^{p^{\prime}}(w)} = 1}}
\nint f(x) S(hw)(x) \sigma \\
&\geq& \nint f(x) S(\mathbf{1}_{[0,1)} w)(x) w([0,1))^{\frac{-1}{p^{\prime}}} \sigma
\end{eqnarray*}
so that
\begin{eqnarray*}
[w]_{A_p} &\lesssim& w([0,1))^{\frac{-1}{p^{\prime}}} \left( \nint_{[0,1)} S(\mathbf{1}_{[0,1)}w)(x)^{p^{\prime}} \sigma \right)^{\frac{1}{p^{\prime}}} \\ &\lesssim& \norme{S(\cdot \sigma)}_{L^p(\sigma) \rightarrow L^p(w)} \\
&\sim& \norme{S}_{L^p(w) \rightarrow L^p(w)}. 
\end{eqnarray*}
As a result, $[w]_{A_p}^{\max \{ 1 , \frac{1}{p-1} \} } \lesssim \norme{S}_{L^p(w) \rightarrow L^p(w)}$. Since $S$ is a positive operator, $S$ extends to a vector-valued operator $\mathcal{S}$ on $L^p_{\ell^r}(w)$ defined by 
\begin{eqnarray*}
\mathcal{S}(\vecc{f})(x) &=& \left( \displaystyle\sum_{j=1}^{\infty} |S(f_j)(x)|^r \right)^{\frac{1}{r}}
\end{eqnarray*}
and $\norme{\mathcal{S}}_{L^p_{\ell^r}(w) \rightarrow L^p(w)} \sim [w]_{A_p}^{\max \{ 1, \frac{1}{p-1} \} }$. 
\end{subsection}
\end{section}
\begin{section}{Proof of Theorem \ref{t.main4}}
\noindent
Let $T$ be as in the statement of Theorem \ref{t.main4}. For each $j$ we apply Lerner's inequality to obtain the following bound:
\begin{eqnarray*}
T(f_j)(x) &\lesssim& Mf_j(x) + \displaystyle\sum_{Q \in \mathcal Q_j} \mathbf{1}_{Q}(x) \omega_{2^{-n-1}}(T(f_j);Q) .
\end{eqnarray*}
For each $j$, we have by Lemma \ref{cz},
\begin{eqnarray*}
\omega_{2^{-n-1}}(T(f_j);Q) &\lesssim& \displaystyle\sum_{m=0}^{\infty} \frac{1}{2^{m \alpha}} \left( \frac{1}{|2^mQ|} \nint_{2^mQ} |f_j(y)| dy \right).
\end{eqnarray*}
Now we make an observation (see \cite{handp},\cite{lerner22}, \cite{handlandp}), for any cube $Q \subset \ral^n$ there is $u$ and $I \in \mathcal{D}^{u}$ such that $Q \subset I$ and $\ell(I) \leq 6 \ell(Q)$. Hence for each $u \in \{ 0, 3^{-1} \}^n$ we may choose a collection of dyadic cubes $\mathcal{Q}_{j, u}$ in $\mathcal D^{u}$ such that
\begin{eqnarray*}
\displaystyle\sum_{Q \in \mathcal Q_j} \mathbf{1}_{Q}(x) \mathbb{E}_{2^mQ}(f_j) &\lesssim& 
\displaystyle\sum_{u \in \{ 0, 3^{-1} \} } \displaystyle\sum_{Q \in \mathcal Q_{j,u}} \mathbf{1}_{Q}(x) \mathbb{E}_{Q}(f_j) \\
&=& \displaystyle\sum_{u \in \{ 0, 3^{-1} \}} P_{j,m,u}(f_j)(x).
\end{eqnarray*}
Define
\begin{eqnarray*}
\vecc{P}_{m,u}(\vecc{f})(x) &=& \left( \displaystyle\sum_{j=1}^{\infty} |P_{j,m,u}(f_j)(x)|^r \right)^{\frac{1}{r}};
\end{eqnarray*}
we have the following bound:
\begin{eqnarray}
\vop{T}(\vecc{f})(x) &\lesssim&  M_r(\vecc{f})(x) + \displaystyle\sum_{m=0}^{\infty} \frac{1}{2^{\alpha m}} \displaystyle\sum_{u \in \{ 0, 3^{-1} \}^n } \nint_{\ral^n}  \vecc{P}_{m,u}(\vecc{f})(x). \label{strtt}
\end{eqnarray}
By Theorem \ref{t.maxv}
\begin{eqnarray*}
\nint_{\ral^n} \vecc{M}_r(\vecc{f})(x)^p w &\lesssim& [w]_{A_p}^{\max \{ \frac{p}{r}, \frac{p}{p-1} \} } \nint_{\ral^n} \norme{\vecc{f}}_{\ell^r}^p w. 
\end{eqnarray*}
For fixed $m$ and $u$, we may apply Theorem \ref{t.3pt}
to obtain
\begin{eqnarray*}
\nint_{\ral^n} \vecc{P}_{m,u}(\vecc{f})(x)^p w &\lesssim& m^4 [w]_{A_p}^{\max \{ p, \frac{p}{p-1} \} } \nint_{\ral^n} \norme{\vecc{f}}_{\ell^r}^p w
\end{eqnarray*}
so that
\begin{eqnarray*}
\displaystyle\sum_{m=0}^{\infty} \frac{1}{2^{\alpha m}} \displaystyle\sum_{u \in \{ 0, 3^{-1} \}^n } \nint_{\ral^n}  \vecc{P}_{m,u}(\vecc{f})(x)^p w &\lesssim&  \left(  \displaystyle\sum_{m=0}^{\infty} \frac{m^4 }{2^{\alpha m}} \right) [w]_{A_p}^{\max \{ p, \frac{p}{p-1} \} } \nint_{\ral^n} \norme{\vecc{f}}_{\ell^r}^p w  \\
&\lesssim& [w]_{A_p}^{\max \{ p, \frac{p}{p-1} \} } \nint_{\ral^n} \norme{\vecc{f}}_{\ell^r}^p w.
\end{eqnarray*}
As a result, from $(\ref{strtt})$ we have
\begin{eqnarray*}
\nint_{\ral^n} \vop{T}(\vecc{f})(x)^p w &\lesssim& [w]_{A_p}^{\max \{ p, \frac{p}{p-1} \} } \nint_{\ral^n} \norme{\vecc{f}}^p w . 
\end{eqnarray*}
\end{section}

\end{document}